\documentclass[12pt, reqno]{amsart}
\usepackage{amsmath, amsthm, amscd, amsfonts, amssymb, graphicx, color}
\usepackage[english]{babel}
\usepackage{cite}
\usepackage[bookmarksnumbered, colorlinks, plainpages]{hyperref}

\newtheorem{theorem}{Theorem}[section]
\newtheorem{lemma}[theorem]{Lemma}

\theoremstyle{definition}

\theoremstyle{remark}

\numberwithin{equation}{section}
\begin{document}
\setcounter{page}{1}

\title[Conditions of coincidence of  central extensions]{Conditions
of coincidence  of central extensions
of von Neumann algebras and algebras of measurable operators}

\author[S.~Albeverio,  K.~K.~Kudaybergenov,
R.~T.~Djumamuratov]{S.~Albeverio$^1$,   K.~K.~Kudaybergenov$^2$
 and R.~T.~Djumamuratov$^3$}

\address{$^{1}$ Institut f\"{u}r Angewandte Mathematik and HCM,
Rheinische Friedrich-Wilhelms-Uni\-versit\"{a}t Bonn, Endenicher
Allee 60,  D-53115 Bonn, Germany}

\email{\textcolor[rgb]{0.00,0.00,0.84}{albeverio@uni-bonn.de}}

\address{$^{2}$ Department of Mathematics, Karakalpak state university\\
Ch. Abdirov 1,  230113, Nukus,    Uzbekistan
\newline}

\email{\textcolor[rgb]{0.00,0.00,0.84}{karim2006@mail.ru}}

\address{$^{3}$ Department of Mathematics, Karakalpak state university\\
Ch. Abdirov 1,  230113, Nukus,    Uzbekistan
\newline}

\email{\textcolor[rgb]{0.00,0.00,0.84}{rauazh@mail.ru}}




\begin{abstract}
Given a von Neumann algebra $M$ we  consider the  central
exten\-sion $E(M)$ of $M.$
We describe class of von Neumann algebras $M$ for which the algebra   $E(M)$ coincides
with the algebra $S(M)$ -- the algebra of all measurable operators
 with respect to   $M,$ and  with $S(M, \tau)$
 -- the algebra of all $\tau$-measurable operators
 with respect to   $M.$
\end{abstract}
 \maketitle

\section{Introduction}

In the series of paper \cite{Alb2, AK1, AK3, AK2}  we have considered
derivations on the algebra $LS(M)$ of locally measurable operators
affiliated with a von Neumann algebra $M,$ and on various
subalgebras of $LS(M).$ A complete description of derivations has
been obtained in the case of  von Neumann algebras of type I and
III.
A comprehensive survey of recent results concerning derivations on
various algebras of unbounded operators affiliated with von
Neumann algebras is presented in \cite{AK2}.
The  general form  of automorphisms  on the algebra $LS(M)$  in the case of
 von Neumann algebras of type I has been obtained in \cite{AK3}. In the proof
 of the main results of the above  papers the crucial role is played by the
central extensions  of von Neumann algebras and also by various topologies
 considered in \cite{AK1}.

 Let  $M$ be an arbitrary von Neumann algebra with the center $Z(M)$
 and let $LS(M)$ denote the algebra of all locally measurable operators with
 respect $M.$ We  consider the set   $E(M)$  of all elements  $x$ from  $LS(M)$
  for which there exists a sequence of
mutually orthogonal central projections  $\{z_i\}_{i\in I}$ in  $M$ with
 $\bigvee\limits_{i\in I}z_i=\textbf{1},$
such that $z_i x\in M$ for all $i\in I.$
It is known  \cite{AK1} that  $E(M)$ is a *-subalgebra in  $LS(M)$ with the center
 $S(Z(M)),$ where   $S(Z(M))$ is  the algebra of all measurable operators
 with respect to   $Z(M),$ moreover,
  $LS(M)=E(M)$ if and only if $M$ does not have
direct summands of type II.
A similar notion (i.e. the algebra $E(\mathcal{A})$) for
arbitrary *-subalgebras $\mathcal{A}\subset LS(M)$  was independently
introduced  by M.A. Muratov and V.I. Chilin \cite{Mur1}.
The algebra  $E(M)$ is called
\textit{the central extension of} $M.$

In section 2  we recall the notions of the algebras $S(M)$ of
measurable operators and $LS(M)$ of locally measurable operators
affiliated with a von Neumann algebra $M.$ We also consider the
\textit{central extension} $E(M)$ of the von Neumann
algebra $M.$

In section 3 we describe the classes of von Neumann algebras $M$ for which the algebra
  $E(M)$ coincides
with the algebras $S(M), S(M, \tau)$ and $M.$

\section{Central extensions of von Neumann algebras}

In this section  we recall the notions of the algebras $S(M)$ of
measurable operators and respectively $LS(M)$ of locally measurable operators
affiliated with a von Neumann algebra $M.$ We also consider the
 \textit{central extension} $E(M)$ of the von Neumann
algebra $M.$

Let  $H$ be a complex Hilbert space and let  $B(H)$ be the algebra
of all bounded linear operators on   $H.$ Consider a von Neumann
algebra $M$  in $B(H)$ with the operator norm $\|\cdot\|_M.$ Denote by
$P(M)$ the lattice of projections in $M.$

A linear subspace  $\mathcal{D}$ in  $H$ is said to be
\emph{affiliated} with  $M$ (denoted as  $\mathcal{D}\eta M$), if
$u(\mathcal{D})\subset \mathcal{D}$ for every unitary  $u$ from
the commutant
$$M'=\{y\in B(H):xy=yx, \,\forall x\in M\}$$ of the von Neumann algebra $M.$

A linear operator  $x$ on  $H$ with the domain  $\mathcal{D}(x)$
is said to be \emph{affiliated} with  $M$ (denoted as  $x\eta M$) if
$\mathcal{D}(x)\eta M$ and $u(x(\xi))=x(u(\xi))$
 for all  $\xi\in
\mathcal{D}(x)$ and for every unitary  $u\in M'.$

A linear subspace $\mathcal{D}$ in $H$ is said to be \emph{strongly
dense} in  $H$ with respect to the von Neumann algebra  $M,$ if

1) $\mathcal{D}\eta M;$

2) there exists a sequence of projections
$\{p_n\}_{n=1}^{\infty}$ in $P(M)$  such that
$p_n\uparrow\textbf{1},$ $p_n(H)\subset \mathcal{D}$ and
$p^{\perp}_n=\textbf{1}-p_n$ is finite in  $M$ for all
$n\in\mathbb{N},$ where $\textbf{1}$ is the identity in $M.$

A closed linear operator  $x$ acting in the Hilbert space $H$ is said to be
\emph{measurable} with respect to the von Neumann algebra  $M,$ if
 $x\eta M$ and $\mathcal{D}(x)$ is strongly dense in  $H.$

 Denote by $S(M)$  the set of all linear operators on $H,$ which are
 measurable with
respect to the von Neumann algebra $M.$ If $x\in S(M),$
$\lambda\in\mathbb{C},$ where $\mathbb{C}$  is the field of
complex numbers, then $\lambda x\in S(M)$  and the operator
$x^\ast,$  adjoint to $x,$  is also measurable with respect to $M$
(see \cite{Seg}). Moreover, if $x, y \in S(M),$ then the operators
$x+y$  and $xy$  are defined on dense subspaces and admit closures
that are called, correspondingly, the strong sum and the strong
product of the operators $x$  and $y,$  and are denoted by
$x\stackrel{.}+y$ and $x \ast y,$
respectively.  It was shown in  \cite{Seg} that
$x\stackrel{.}+y$ and $x \ast y$ belong to $S(M)$ and
these algebraic operations make $S(M)$ a $\ast$-algebra with the
identity $\textbf{1}$  over the field $\mathbb{C}.$ Here, $M$ is a
$\ast$-subalgebra of $S(M).$ In what follows, the strong sum and
the strong product of operators $x$ and $y$  will be denoted in
the same way as the usual operations, by $x+y$  and $x y,$
respectively.

A closed linear operator $x$ in  $H$  is said to be \emph{locally
measurable} with respect to the von Neumann algebra $M,$ if $x\eta
M$ and there exists a sequence $\{z_n\}_{n=1}^{\infty}$ of central
projections in $M$ such that $z_n\uparrow\textbf{1}$ and $z_nx \in
S(M)$ for all $n\in\mathbb{N}$ (see \cite{Yea}).

Denote by $LS(M)$ the set of all linear operators that are locally
measurable with respect to $M.$ It was proved in \cite{Yea}  that
$LS(M)$ is a $\ast$-algebra over the field $\mathbb{C}$ with
identity $\textbf{1},$ the operations of strong addition, strong
multiplication, and passing to the adjoint. In such a case, $S(M)$
is a $\ast$-subalgebra in $LS(M).$ In the case where $M$ is a
finite von Neumann algebra or a factor, the algebras $S(M)$ and
$LS(M)$ coincide. This is not true in the general case. In
\cite{Mur2}  the class of von Neumann algebras $M$ has been described for
which the algebras  $LS(M)$ and  $S(M)$ coincide.

 Let   $\tau$ be a faithful normal semi-finite trace on $M.$ We recall that
  a closed linear operator
  $x$ is said to be  $\tau$\textit{-measurable} with respect to the von Neumann algebra
   $M,$ if  $x\eta M$ and   $\mathcal{D}(x)$ is
  $\tau$-dense in  $H,$ i.e. $\mathcal{D}(x)\eta M$ and given   $\varepsilon>0$
  there exists a projection   $p\in M$ such that   $p(H)\subset\mathcal{D}(x)$
  and $\tau(p^{\perp})<\varepsilon.$
   Denote by  $S(M,\tau)$ the set of all   $\tau$-measurable operators with respect to  $M$
   (see \cite{Nel}).

    It is well-known that  $S(M, \tau)$ is a $\ast$-subalgebras in $LS(M)$ (see \cite{Nel}).

    Consider the topology  $t_{\tau}$ of convergence in measure or \textit{measure topology}
    on $S(M, \tau),$ which is defined by
 the following neighborhoods of zero:
\begin{center}
$
V(\varepsilon, \delta)=\{x\in S(M, \tau): \exists\,
 e\in P(M), \tau(e^{\perp})\leq\delta, xe\in
M,  \|xe\|_{M}\leq\varepsilon\},
$
\end{center}
  where $\varepsilon, \delta$
are positive numbers, and $\|\cdot\|_{M}$ denotes the operator norm on
$M$.

 It is well-known \cite{Nel}) that $S(M, \tau)$ equipped with the measure topology is a
complete metrizable topological $\ast$-algebra.

Let $(\Omega,\Sigma,\mu)$  be a measure space and
suppose
 that the measure $\mu$ has the  direct sum property, i.e. there is a family
 $\{\Omega_{i}\}_{i\in
J}\subset\Sigma,$ $0<\mu(\Omega_{i})<\infty,\,i\in J,$ such that
for any $A\in\Sigma,$ $\mu(A)<\infty,$ there exist a countable
subset
$J_{0 }\subset J$ and a set  $B$ with zero measure such
that  $A=\bigcup\limits_{i\in J_{0}}(A\cap
\Omega_{i})\cup B.$

It is well-known (see e.g.~\cite{Seg}) that every commutative von
Neumann algebra
 $M$
is $\ast$-isomorphic to the algebra $L^{\infty}(\Omega, \Sigma, \mu)$
of all (equivalence classes of) complex essentially bounded
measurable functions on  $(\Omega, \Sigma, \mu)$ and in this case
$LS(M)=S(M)\cong L^{0}(\Omega, \Sigma, \mu),$ where $L^{0}(\Omega,
\Sigma, \mu)$ the algebra of all (equivalence classes of) complex
measurable functions on $(\Omega, \Sigma, \mu).$

Further we consider the algebra  $S(Z(M))$  of operators which are measurable
with respect to the  center $Z(M)$ of the von Neumann algebra $M.$
Since  $Z(M)$ is an abelian von Neumann algebra  it
 is $\ast$-isomorphic to $ L^{\infty}(\Omega, \Sigma, \mu)$
   for an appropriate measure space $(\Omega, \Sigma, \mu)$.
   Therefore the algebra  $S(Z(M))$ coincides
   with $Z(LS(M))$ and  can be
 identified with the algebra $ L^{0}(\Omega, \Sigma, \mu).$

The basis of neighborhoods of zero in the topology of convergence locally in measure
    on $L^0(\Omega,\Sigma, \mu)$ consists of the sets
$$W(A,\varepsilon,\delta)
=\{f\in L^0(\Omega,\Sigma, \mu):\exists B\in \Sigma, \, B\subseteq A, \,
\mu(A\setminus B)\leq \delta, $$
$$ f\cdot \chi_B \in L^{\infty}(\Omega,\Sigma, \mu),\,
||f\cdot \chi_B||_{L^{\infty}(\Omega,\Sigma, \mu)}\leq \varepsilon\},$$
where $\varepsilon, \delta>0, \, A\in \Sigma, \, \mu(A)<+\infty,$
and $\chi_B$ is the characteric
function of the set $B\in \Sigma.$

Let us recall the definition of the dimension functions $d$ on the lattice
$P(M)$ of projection from $M$ (see \cite{Mur}, \cite{Seg}).

Let   $L_+$  denote the set of all measurable  functions  $f:
(\Omega,\Sigma, \mu)\rightarrow [0,{\infty}]$ (modulo functions
equal to zero $\mu$-almost everywhere).

 Let  $M$ be an arbitrary von Neumann algebra with the center
  \linebreak $Z(M)\equiv L^\infty(\Omega,\Sigma, \mu).$
Then there exists a map  $d:P(M)\rightarrow L_{+}$ with the
following properties:

(i) $d(e)$ is a finite function if only if the projection  $e$ is finite;

(ii) $d(e+q)=d(e)+d(q)$  for $p, q \in P(M),$ $eq=0;$

(iii) $d(uu^*)=d(u^*u)$ for every  partial isometry  $u\in M;$

(iv) $d(ze)=zd(e)$ for all $z\in P(Z(M)), \,\, e\in P(M);$

(v) if  $\{e_{\alpha}\}_{\alpha \in J}, \,\,\, e\in P(M) $ and $e_{\alpha}\uparrow e,$ then
$d(e)=\sup \limits_{\alpha \in J}d(e_{\alpha}).$

This map  $d:P(M)\rightarrow L_+,$ is a called the \emph{dimension functions} on  $P(M).$

The basis of neighborhoods of zero in the topology $t(M)$ of
\emph{convergence locally in measure}  on $LS(M)$ consists (in the
above notations) of the following sets
\begin{align*}
V(A,\varepsilon,\delta)=\{x\in LS(M):\exists p\in P(M), \, \exists
z\in P(Z(M)),  \, xp \in M, \\
 ||xp||_{M}\leq \varepsilon, \,\, z^{\bot}\in
W(A,\varepsilon,\delta), \,\, d(zp^{\bot})\leq \varepsilon z\},
\end{align*} where
 $\varepsilon, \delta>0, \, A\in \Sigma, \, \mu(A)<+\infty$ (see
 \cite{Yea}).

The topology  $t(M)$ is  metrizable if and only if the center  $Z(M)$
is   $\sigma$-finite (see \cite{Mur}).

Given an arbitrary  family  $\{z_i\}_{i\in I}$ of mutually orthogonal
central projections in $M$ with $\bigvee\limits_{i\in
I}z_i=\textbf{1}$ and a  family of elements $\{x_i\}_{i\in I}$ in
$LS(M)$ there exists a unique element $x\in LS(M)$ such that $z_i
x=z_i x_i$ for all $i\in I.$ This element is denoted by
$x=\sum\limits_{i\in I}z_i x_i.$

We  denote by  $E(M)$  the set of all elements  $x$ from  $LS(M)$
 for which there exists a sequence of
mutually orthogonal central projections  $\{z_i\}_{i\in I}$ in
 $M$ with $\bigvee\limits_{i\in I}z_i=\textbf{1},$
such that $z_i x\in M$ for all $i\in I,$ i.e.
 $$E(M)=\{x\in LS(M): \exists
  z_i\in P(Z(M)), z_iz_j=0, i\neq j, \bigvee\limits_{i\in I}z_i=\textbf{1},
 z_i x\in M, i\in I\},$$
where $Z(M)$ is the center of $M.$

It is known  \cite{AK1} that  $E(M)$ is  *-subalgebras in  $LS(M)$ with the center
 $S(Z(M)),$ where   $S(Z(M))$ is  the algebra of all measurable operators
 with respect to   $Z(M),$ moreover,
  $LS(M)=E(M)$ if and only if $M$ does not have
direct summands of type II.

A similar notion (i.e. the algebra $E(\mathcal{A})$) for
arbitrary *-subalgebras $\mathcal{A}\subset LS(M)$  was independently
introduced recently by M.A. Muratov and V.I. Chilin \cite{Mur1}.
The algebra  $E(M)$ is called
\textit{the central extension of} $M.$

It is known \cite{AK1},
\cite{Mur1} that an element
$x\in  LS(M)$ belongs to $E(M)$ if and only if there exists
 $f\in S(Z(M))$ such that    $|x|\leq f.$
Therefore for each
 $x\in E(M)$ one can define the following vector-valued norm
$$
 ||x||=\inf\{f\in S(Z(M)): |x|\leq f\}
$$
and this norm satisfies the following conditions:

$1) \|x\|\geq 0; \|x\|=0\Longleftrightarrow x=0;$

$2)  \|f x\|=|f|\|x\|;$

$ 3)  \|x+y\|\leq\|x\|+\|y\|;$

$4) ||x y||\leq ||x||||y||;$

$5) ||xx^{\ast}||=||x||^2$
\newline
    for all $x,y\in E(M), f\in S(Z(M)).$

 Let  $M$ be an arbitrary von Neumann algebra with the
 center $Z(M)\equiv L^\infty(\Omega,\Sigma, \mu).$
On the space  $E(M)$ we consider the following sets:
$$
O(A, \varepsilon, \delta)=\left\{x \in E(M): ||x||\in W(A, \varepsilon, \delta)\right\},
$$
where  $\varepsilon ,\delta > 0,\,\,\,A \in \sum ,\,\,\,\mu \left( A \right) <
+ \infty $.

The system  of sets
\begin{equation}
 \label{base2}
\{x+O(A, \varepsilon, \delta)\},
\end{equation}
where $x\in  E(M), \varepsilon>0, \delta>0,
A\in \Sigma, \mu(A)<\infty$, defines on  $E(M)$ a Hausdorff vector topology $t_c(M),$
for which the sets (\ref{base2}) form the base of neighborhoods of the element $x\in E(M).$
Moreover  in this topology the   involution is continuous and the
 multiplication  is jointly continuous,
i.e. $(E(M), t_c(M))$ is a topological $\ast$-algebra.
It is known  \cite[Proposition 3.2]{AKJ} that  $(E(M), t_c(M))$ is
 a complete topological $\ast$-algebra and
 $M$ is a  $t_c(M)$-dense in  $E(M).$

\section{Conditions of coincidence of central extensions
of von Neumann algebras and algebras of measurable operators}

In this section  we describe class of von Neumann algebras $M$ for
which the algebra   $E(M)$ coincide
with algebra $S(M), S(M, \tau)$ and $M.$

It should be noted that  \cite[Proposition 1.1]{AK1} and \cite[Theorem 3.1]{AKJ}
 imply the following result.

\begin{theorem}\label{TE}
The following conditions on a given von Neumann algebra $M$ are equivalent:

(1) $E(M)=LS(M);$

(2)  $M$ does not have
direct summands of type II.

In this case  the topologies $t_c(M)$ and $t(M)$
coincide.
\end{theorem}

Now  we   describe a class of von Neumann algebras $M$ for which the algebras  $E(M)$ and $S(M)$
coincide.

\begin{theorem}\label{TA}
The following conditions on a given von Neumann algebra $M$ are equivalent:

(1) $E(M)=S(M);$

(2)  $M=\bigoplus\limits_{k=0}^{n}M_k,$ where
$M_0$ be a von Neumann algebra of type  $I_{fin},$ $M_k$ be a factors of type   $I_\infty$
or III, $k=\overline{1, n}.$

In this case  the topologies $t_c(M)$ and $t(M)$
coincide.
\end{theorem}

For the proof of Theorem \ref{TA} we
 need following result.

\begin{lemma}\label{LA}
If   $M$ be a  von Neumann algebra of type II with a faithful normal semifinite trace
$\tau,$ then
 $E(M)\neq S(M, \tau)$ and  $E(M)\neq S(M).$
 \end{lemma}

\begin{proof}
Suppose that $M$ is a type II von Neumann
algebra. First assume that $M$ is of type II$_1$  and admits a
faithful normal tracial state $\tau$ on $M.$ Without loss generality we assume that
$\tau(\textbf{1})=1.$ Let  $\Phi$ be the
canonical center-valued trace on $M.$
Since   $M$ is of type  II,   there exists a projection
$p_1\in M$ such that  $$p_1\sim \textbf{1}-p_1.$$ Then
$\Phi(p_1)=\Phi(p_1^{\perp}).$ From
$\Phi(p_1)+\Phi(p_1^{\perp})=\Phi(\textbf{1})=\textbf{1}$ it follows that
$$\Phi(p_1)=\Phi(p_1^{\perp})=\frac{\textstyle 1}{\textstyle 2}\textbf{1}.$$

Suppose that there exist  mutually orthogonal projections
$p_1,\,p_2,\cdots,p_n$ in $M$ such that
$$\Phi(p_k)=\frac{\textstyle 1}{\textstyle 2^{k}}\textbf{1},\,k=\overline{1, n}.$$
Set  $e_n=\sum\limits_{k=1}^{n}p_k.$ Then
$\Phi(e_n^{\perp})=\frac{\textstyle 1}{\textstyle
2^{n}}\textbf{1}.$ Take a projection  $p_{n+1}<e_n^{\perp}$ such that
$$p_{n+1}\sim e_n^{\perp}-p_{n+1}.$$
Then
$$\Phi(p_{n+1})=\frac{\textstyle 1}{\textstyle
2^{n+1}}.$$

Hence there exists a sequence a mutually orthogonal
projections $\{p_n\}_{n\in\mathbb{N}}$ in $M$ such that
$$\Phi(p_n)=\frac{\textstyle
1}{\textstyle 2^{n}}\textbf{1},\,n\in\mathbb{N}.$$
Note that  $\tau(p_n)=\frac{\textstyle 1}{\textstyle 2^{n}}.$
Indeed
$$\tau(p_n)=\tau(\Phi(p_n))=\tau\left(\frac{\textstyle
1}{\textstyle 2^{n}}\textbf{1}\right)=\frac{\textstyle 1}{\textstyle
2^{n}}.$$

Since
$$\sum\limits_{n=1}^{\infty}n\tau(p_n)=
\sum\limits_{n=1}^{\infty}\frac{n}{2^{n}}<+\infty$$
it follows that  the series
$$\sum\limits_{n=1}^{\infty}np_n$$ converges in measure in  $S(M, \tau).$
Therefore
there exists  $x=\sum\limits_{n=1}^{\infty}np_n\in S(M, \tau).$

Let us  show that $x\in S(M, \tau)\setminus E(M).$ Suppose that  $zx\in M,$
where  $z$ is a nonzero central
projection. Since any  $p_n$ is a faithful projection we have that
$z p_n\neq 0$ for all $n.$ Thus
$$||zx||_{M}=1||zx||_{M}1=||p_n||_{M}\cdot||zx||_{M}\cdot||p_n||_{M}\geq||zp_nxp_n||_{M}=
||zp_n n||_{M}=n,$$
i.e.
$$||zx||_{M}\geq n$$
for all $n\in\mathbb{N}.$ From this  contradiction it follows that
 $x\in S(M, \tau)\setminus E(M).$

  For a general type II von Neumann
 algebra $M$ take a non zero finite projection $e\in M$ and
 consider the finite type II von Neumann algebra $eMe$ which
 admits a separating family of normal tracial states. Now if $f\in
 eMe$ is the support projection of some tracial state $\tau$ on
 $eMe$ then $fMf$ is a type II$_1$ von Neumann algebra with a
 faithful normal tracial state. Hence as above one can
 construct an element $x\in S(M, \tau)\setminus E(M),$ moreover
 $x\in S(M)\setminus E(M).$
 The proof is complete. \end{proof}

\textit{The proof of the theorem \ref{TA}.}  (1) $\Rightarrow$ (2).
If the algebra $M$ has a direct summand of type  II, then by  Lemma \ref{LA} we have that
 $E(M)\neq S(M).$ Hence if
 $E(M)=S(M),$
then there exist mutually orthogonal central projections
$z_0, z_1, z_2$  with  $z_0+z_1+z_2=\textbf{1}$ such that
$z_0M$ is a type  I$_{fin},$ $z_1M$ is a type
I$_\infty$ and $z_2M$ is a type  III.

Suppose that $zZ(M)$ is infinite dimensional,
where  $Z(M)$ is the  center  $M$ and $z=z_1+z_2.$ Then there exists a sequence
of nonzero mutually orthogonal projections  $\{p_n\}_{n=1}^\infty$
in  $zZ(M).$ Put
\begin{equation}
\label{EA}
x=\sum\limits_{n=1}^\infty n p_n.
\end{equation}
Then  $0\leq x\in E(M)$ and  $e_n(x)=\sum\limits_{k=1}^n  p_k,$ where
$e_n(x)$ is  a spectral projection   $x$ corresponding
to the interval  $[0, n].$
Since  $z M$ is a properly infinite algebra, then  $e_n(x)^\perp=\sum\limits_{k=n+1}^\infty  p_k$ is
an infinite projection for all  $n\in\mathbb{N}.$  This means that
$x\notin S(M),$  i.e. $E(M)\neq S(M).$ This is a contradiction with $E(M)=S(M).$
Therefore
$zZ(M)$ is a finite dimensional algebra. Thus
$$
zM=\bigoplus\limits_{k=1}^n M_k,$$
where $M_k$ is a factor of type  I$_\infty$ or
III, $k=\overline{1, n}.$ Put  $M_0=z_0 M.$ Then
$$
M=z_0 M\oplus z M=\bigoplus\limits_{k=0}^n M_k,$$
where  $M_0$ is a type
I$_{fin},$
$M_k$ is a factor of type  I$_\infty$ or
III, $k=\overline{1, n}.$

(2) $\Rightarrow$ (1)
Let
$$
M=\bigoplus\limits_{k=0}^n M_k,$$
where  $M_0$ be a type
I$_{fin},$
$M_k$ be a factor of type  I$_\infty$ or
III, $k=\overline{1, n}.$
 Since  $M_0$ is  a type
I$_{fin},$ then from  \cite[Proposition  1.1]{AK1}
it follows that
$$
E(M_0)=LS(M_0)=S(M_0).
$$ Since  $M_k$ are   factors of type  I$_\infty$ or
III, $k=\overline{1, n},$
then by  \cite[Theorem  1]{Mur2}  we obtain that
$$
S(M_k)=M_k=E(M_k),
$$
for all $k=\overline{1, n}.$ Hence
$$
E(M)=\bigoplus\limits_{k=0}^n E(M_k)=
\bigoplus\limits_{k=0}^n S(M_k)=S(M),
$$
i.e.
$E(M)=S(M).$

Now let
$M=\bigoplus\limits_{k=0}^{n}M_k,$ where
$M_0$ is  a von Neumann algebra of type  I$_{fin},$ $M_k$
are  factors of type   I$_\infty$
or III, $k=\overline{1, n}.$ Then $LS(M)=S(M)=E(M)$ and
$M$ does not have
direct summands of type II. Therefore from Theorem \ref{TE} we obtain that the topologies
$t(M)$ and $t_c(M)$ coincide.
The proof is complete. $\Box$

We now  describe a class of von Neumann algebras $M$ for which the algebras
$E(M)$ and  $S(M, \tau)$
coincide.

\begin{theorem}
Let  $M$ be a von Neumann algebra with a faithful normal semifinite
trace
$\tau.$
The following conditions  are equivalent:

(1) $E(M)=S(M, \tau);$

(2)  $M=\bigoplus\limits_{k=0}^{n}M_k,$ where
$M_0$ be a type  $I_{fin}$ algebra,  $M_k$ be a factors of type
$I_\infty,$ $k=\overline{1, n},$
and restriction of trace  $\tau$ on  $M_0$ is  finite.

In this case the topologies $t_c(M)$ and $t_\tau$ coincide.
\end{theorem}

\begin{proof}   (1) $\Rightarrow$ (2).
If  $M$ has a direct summand of type  II,
then by Lemma \ref{LA} we obtain that  $E(M)\neq S(M, \tau).$ Hence if
 $E(M)=S(M, \tau),$
 then there exist
orthogonal central projections
$z_0, z_1$ with  $z_0+z_1=\textbf{1}$ such that
$z_0M$ is a type  I$_{fin},$ $z_1M$ is a type
I$_\infty.$

If we assume that  $z_1Z(M)$ is  infinite dimensional then the
element  $x\in E(M)$ defined similarly as  in  (\ref{EA})
we have that $x\notin S(M, \tau),$  i.e.  $E(M)\neq S(M, \tau).$ Hence,
$z_1Z(M)$ is finite dimensional. Thus
$$
zM=\bigoplus\limits_{k=1}^n M_k,$$
where  $M_k$ is of  type  I$_\infty,$ $k=\overline{1, n}.$ Put $M_0=z_0 M.$ Then
$$
M=z_0 M\oplus z_1 M=\bigoplus\limits_{k=0}^n M_k,$$
where  $M_0$ is of  type
I$_{fin},$
$M_k$ are  factors of type  I$_\infty,$ $k=\overline{1, n}.$
Now by Theorem~\ref{TA} we have that
$E(M_0)=S(M_0).$ At the same time by conditions of theorem it follows that
  $E(M_0)=S(M_0, \tau_0),$ where $\tau_0$ is the restriction of $\tau$ on
  $M_0.$ Thus  $S(M_0)=S(M_0, \tau_0).$
Therefore by  \cite[Proposition 9]{Mur2} the restriction of trace  $\tau$ on $M_0$ is finite.

(2) $\Rightarrow$ (1).
Let
$$
M=\bigoplus\limits_{k=0}^n M_k,$$
where  $M_0$ is  an algebra of  type
I$_{fin},$
$M_k$ are  factors of type  I$_\infty,$  $k=\overline{1, n}$
 and the restriction $\tau$ on  $M_0$ be a finite.
 Let    $\tau_0$ be the  restriction  of  $\tau$ on  $M_0.$ Since  $M_0$ is of
  type
 I$_{fin}$ then
 by  \cite[Proposition  1.1]{AK1}
it follows that
$$
E(M_0)=LS(M_0)=S(M_0).
$$
Now  since the  trace  $\tau_0$  is a finite then
$S(M_0)=S(M, \tau_0).$ Thus $E(M_0)=S(M_0, \tau_0).$  Since  $M_k$ is  a
factor of type  I$_\infty,$ $k=\overline{1, n},$
then from  \cite[Theorem~2.2.9]{Mur} we obtain that
$S(M_k, \tau_k)=M_k=E(M_k),$  where $\tau_k$ is the restriction of $\tau$ on
  $M_k,$
$k=\overline{1, n}.$ Therefore
$$
E(M)=\bigoplus\limits_{k=0}^n E(M_k)=
\bigoplus\limits_{k=0}^n S(M_k, \tau_k)=S(M, \tau),
$$
i.e.
$E(M)=S(M, \tau).$

Now let $M=\bigoplus\limits_{k=0}^{n}M_k,$ where
$M_0$ is  an algebra of  type  I$_{fin},$ $M_k$ are  factors of type
  I$_\infty,$ $k=\overline{1, n},$
and the restriction $\tau_0$ of the trace  $\tau$ on  $M_0$ is  finite.
Since the restriction of the trace  $\tau$ on  $M_0$ is  finite then
$E(M_0)=S(M_0)=S(M_0, \tau_0)$ and the restrictions of
the topologies $t_c(M)$ and $t_{\tau}$ on
$E(M_0)$ coincide. Further since $M_k$ are factors of type
 I$_\infty,$ $k=\overline{1, n},$ then
$$
E\left(\bigoplus\limits_{k=1}^{n}M_k\right)=
S\left(\bigoplus\limits_{k=1}^{n}M_k, \tau|_{\bigoplus\limits_{k=1}^{n}M_k}\right)=
\bigoplus\limits_{k=1}^{n}M_k
$$ and the restrictions  of the topologies $t_c(M)$ and $t_\tau$ on
$\bigoplus\limits_{k=1}^{n}M_k$ coincide
with the uniform topology on $\bigoplus\limits_{k=1}^{n}M_k.$
Therefore  the topologies $t_c(M)$ and $t_\tau$ coincide.
The proof is complete.
\end{proof}

Finally we   describe a
class of von Neumann algebras $M$ for which the algebras  $E(M)$ and  $M$
coincide.

\begin{theorem}
Let  $M$ be a von Neumann algebra.
The following conditions  are equivalent:

 (1) $E(M)=M;$

(2)  $M=\bigoplus\limits_{k=1}^{n}M_k,$ where
$M_k$ are  von Neumann factors for all   $k=\overline{1, n}.$

In this case the topology $t_c(M)$ coincides with uniform topology.
\end{theorem}

\begin{proof}
Let  $E(M)=M.$ Suppose that  $Z(M)$ is  infinite dimensional.
Then there exists a sequence of nonzero central orthogonal projections  $\{p_n\}_{n=1}^\infty$
in  $M.$ Put $$
x=\sum\limits_{n=1}^\infty n p_n.
$$
Then  $x\in E(M)\setminus M,$ this is contradiction. This means that
$Z(M)$ is  finite dimensional. Thus
$$
M=\bigoplus\limits_{k=1}^n M_k,$$
where  $M_k$ are von Neumann factors for all     $k=\overline{1, n}.$

(2) $\Rightarrow$ (1).
Let
$$
M=\bigoplus\limits_{k=1}^n M_k,$$
where  $M_k$ are  von Neumann factors for all     $k=\overline{1, n}.$
Then by the definition of central extensions it follows that
$E(M_k)=M_k.$  Therefore
$$
E(M)=\bigoplus\limits_{k=1}^n E(M_k)=
\bigoplus\limits_{k=1}^n M_k=M,
$$
i.e.
$E(M)=M.$

 Now let $M=\bigoplus\limits_{k=1}^{n}M_k,$ where
$M_k$ are  von Neumann factors for all   $k=\overline{1, n}.$
Since $M_k$ are  von Neumann factors for all   $k=\overline{1, n}$
then the restriction of the topology
$t_c(M)$ on $M_k$ coincides with the uniform topology. Therefore
the topology $t_c(M)$ coincides with the uniform topology.
The proof is complete.
\end{proof}

\section*{Acknowledgments}
 The second  named author would like to acknowledge the
  hospitality of the "Institut fur Angewandte Mathematik", Universitat Bonn (Germany).
 This work is supported in part by
the DFG AL 214/36-1 project (Germany). This work is supported in
part by the German Academic Exchange Service -- DAAD.


\begin{thebibliography}{13}



\bibitem{Alb2} Albeverio~S.,  Ayupov~Sh.~A., Kudaybergenov~K.~K.,
 Structure of derivations on
various algebras of measurable operators for type I von Neumann
algebras, J. Func. Anal.,   256 (2009)  2917--2943.


\bibitem{AK3} Albeverio~S.,  Ayupov~Sh.~A., Kudaybergenov~K.~K.,
Djumamuratov~R.~T.,
Automorphisms
  of central extensions  of type I von Neumann algebras, http://arxiv.org/abs/1104.4698.



\bibitem{AK1} Ayupov~Sh.~A., Kudaybergenov~K.~K.,
Additive derivations on algebras of measurable operators, ICTP,
Preprint,  No  IC/2009/059,   -- Trieste, 2009. -- 16 p. (accepted
in \textit{Journal of operator theory}).

\bibitem{AK2}{Ayupov~Sh.~A., Kudaybergenov~K.~K.,} Derivations on
algebras of measurable operators, Infin. Dimens. Anal.
Quantum Probab. Relat. Top.,  13 (2010)   305--337.


\bibitem{AKJ} Ayupov~Sh.~A., Kudaybergenov~K.~K.,
Djumamuratov~R.~T.,
Topologies on central extensions  of von Neumann algebras, http://arxiv.org/abs/1107.5153.



\bibitem{Mur}  Muratov~M.~A.,   Chilin~V.~I.,  Algebras of measurable and locally measurable
 operators,
 Institute of Mathematics Ukrainian Academy of Sciences,  Kiev, 2007.




\bibitem{Mur2}
Muratov~M.~A., Chilin~V.~I., $\ast$-Algebras of unbounded
operators affiliated with a von Neumann algebra, J. Math. Sciences,
140 (2007), 445--451.





\bibitem{Mur1}  Muratov~M.~A.,   Chilin~V.~I.,  Central extensions of *-algebras of
measurable operators, Doklady AN Ukraine,  2009, 24--28.




\bibitem{Nel} Nelson~E., Notes on non-commutative integration, J. Funct.
Anal. 15 (1974)  103--116.


\bibitem{Seg}  Segal~I.,
A non-commutative extension of abstract integration, Ann.
Math., 57 (1953)   401--457.




  \bibitem{Yea} Yeadon~F.~J.,
 Convergence of measurable operators.  Proc. Camb. Phil.
       Soc., 74 (1973)   257--268.




 \end{thebibliography}
\end{document}